\documentclass{amsart}
\usepackage[numbers]{natbib}
\usepackage{thmtools}
\usepackage[breaklinks]{hyperref}
\usepackage{srcltx}

\newcommand{\norm}[2][y]{\if#1y\left\fi\lVert#2\if#1y\right\fi\rVert} 
\newcommand{\abs}[2][y]{\if#1y\left\fi\lvert#2\if#1y\right\fi\rvert} 
\newcommand{\Real}{\mathbb{R}}

\renewcommand{\natural}{\mathbb{N}} 
\newcommand{\LL}{\mathrm{L}}
\newcommand{\HH}{\mathrm{H}}
\newcommand{\CC}{\mathrm{C}}
\newcommand{\dd}[1][y]{\if#1y\,\fi{\mathrm d}} 
\renewcommand{\div}{\mathrm{div}} 
\newcommand{\poly}{\Pi}
\DeclareMathOperator{\proj}{proj}

\declaretheorem[name=Theorem,numberwithin=section]{thm}
\declaretheorem[name=Lemma,sibling=thm]{lemma}
\declaretheorem[name=Proposition,sibling=thm]{proposition}
\declaretheorem[name=Corollary,sibling=thm]{corollary}
\declaretheorem[name=Remark,sibling=thm,style=remark]{remark}

\begin{document}
\title[Orthogonal polynomial projection in the unit ball]{Orthogonal polynomial projection error measured in Sobolev norms in the unit ball}
\thanks{The author was supported by the CONICYT grant FONDECYT-1130923 and by the Universidad de Concepci\'on VRID-Enlace grant 215.013.040-1.0.}
\author{Leonardo E. Figueroa}
\email{lfiguero@ing-mat.udec.cl}
\address{CI\textsuperscript{2}MA and Departamento de Ingenier\'ia Matem\'atica, Universidad de Concepci\'on, Casilla 160-C, Concepci\'on, Chile}
\begin{abstract}
We study approximation properties of weighted $\LL^2$-orthogonal projectors onto spaces of polynomials of bounded degree in the Euclidean unit ball, where the weight is of the generalized Gegenbauer form $x \mapsto (1-\norm{x}^2)^\alpha$, $\alpha > -1$.
Said properties are measured in Sobolev-type norms in which the same weighted $\LL^2$ norm is used to control all the involved weak derivatives.
The method of proof does not rely on any particular basis of orthogonal polynomials, which allows for a short, streamlined and dimension-independent exposition.
\end{abstract}
\maketitle

\section{Introduction}

It has been known since the early eighties \cite{CQ:1982} that the orthogonal projector $S^0_N$ mapping $\LL^2(-1,1)$ onto the space of univariate polynomials of degree less than or equal to $N$ (equivalently, $S^0_N$ is the operation consisting in truncating the Fourier--Legendre series of its argument at degree $N$) satisfies the bound
\begin{equation}\label{LegendreBound}
(\forall\,u\in\HH^l(-1,1)) \quad \norm{u - S^0_N(u)}_{\HH^1(-1,1)} \leq C N^{3/2-l} \norm{u}_{\HH^l(-1,1)},
\end{equation}
where $C > 0$ depends only on $l$ and $\HH^1(-1,1)$ and $\HH^l(-1,1)$ denote standard Sobolev spaces (see \cite[Ch.~5]{CHQZ-I} for a detailed proof of \eqref{LegendreBound} and its Chebyshev weight and periodic unweighted analogues and \cite{Guo:2000a} for its general Gegenbauer weight analogue).
Recently \cite{Figueroa:arXiv2015} this result was extended to the unit disk for Gegenbauer-type weights.

The purpose of this work is proving a weighted analogue of \eqref{LegendreBound} in the case of the unit ball of any dimension; in order to state it, we introduce now the minimal necessary notation.
Let $B^d$ be the unit ball of $\Real^d$, $\alpha > -1$ and let the weight function $W_\alpha \colon B^d \to \Real$ be defined by $W_\alpha(x) = (1-\norm{x}^2)^\alpha$ with $\norm{\cdot}$ being the Euclidean norm.
We denote by $\LL^2_\alpha$ the weighted Lebesgue space $W_\alpha^{-1/2} \LL^2(B^d)$, whose natural squared norm is $\norm{u}_\alpha^2 := \int_{B^d} \abs{u}^2 W_\alpha$, and, given an integer $l \geq 0$, by $\HH^l_\alpha$ the weighted Sobolev space whose squared norm is $\norm{u}_{\HH^l_\alpha}^2 := \sum_{k=0}^l \norm{\nabla_k u}_\alpha^2$.
Let $S^\alpha_N$ be the orthogonal projector mapping $\LL^2_\alpha$ onto $\poly^d_N$, where $\poly^d_N$ is the space of $d$-variate polynomials of degree less than or equal to $N$.
Our main result is

\begin{thm}\label{thm:lossy}
For all integers $1 \leq r \leq l$ there exists $C = C(d,\alpha,l,r) > 0$ such that
\begin{equation}\label{lossyPreview}
(\forall\,u\in\HH^l_\alpha) \quad \norm{u - S^\alpha_N(u)}_{\HH^r_\alpha} \leq C N^{-1/2 + 2r - l} \norm{u}_{\HH^l_\alpha}.
\end{equation}
\end{thm}

There are two application domains of our main result that we are aware of.
One lies in the analysis of polynomial interpolation operators (cf.\ \cite{CQ:1982} and \cite[Ch.~5]{CHQZ-I}), themselves important in the analysis of spectral methods.
The other, which is the one that led us into this pursuit in the first place, lies in the characterization of approximability spaces relevant to the analysis of nonlinear iterative methods for the numerical solution of high-dimensional PDE; we remit the interested reader to \cite[Ch.~4]{Figueroa} where the one-dimensional case of \autoref{thm:lossy} and the fact that the $S\sp{\alpha}_N$ projectors tensorize in a very straightforward way are exploited for such task.

As noted ever since \cite{CQ:1982}, \autoref{thm:lossy} compares unfavorably with the situation for trigonometric polynomials in unweighted periodic Sobolev spaces, where the power on $N$ is simply $r-l$.
The origin of this difference in behavior is that in the trigonometric case differentiation and projection commute, something which is impossible in the algebraic case \cite[\S~2.3.2]{CHQZ-I}.

We emphasize that the case $r = 0$ is explicitly excluded from consideration in \autoref{thm:lossy}, for in such a case the provably optimal power on $N$ is $-l$ (cf.\ \autoref{lem:L2-projection-error} below), outside the pattern set in \eqref{lossyPreview}.
We also note that if $2r \geq l+1/2$ in \eqref{lossyPreview}, $S^\alpha_N(u)$ need not converge to $u$ in $\HH^r_\alpha$ as $N$ tends to infinity.
We further remark that \autoref{thm:lossy} is not a best or quasi-best approximation result (for those see \cite[Ch.~5]{CHQZ-I}, \cite{Guo:2000a}, \cite[\S~4]{LiXu:2014} and \cite[\S~5]{DaiXu:2011}), because in general the orthogonal projection of $\HH^r_\alpha$ onto $\poly^d_N$ need not coincide with the restriction of $S^\alpha_N$ to $\HH^r_\alpha$.

In every proof of a particular instance of \autoref{thm:lossy} that we are aware of, an important role was played by spectral differentiation formulas, which connect the orthogonal expansion coefficients of a function and one of its derivatives; e.g., \cite[Eq.~(2.3.18)]{CHQZ-I}
\begin{equation*}
(\forall\,k \in \{0, 1, 2, \dotsc \}) \quad \hat u\sp{(1)}_k = (2k+1) \sum_{q=0}^\infty \hat u_{k+1+2q},
\end{equation*}
where $u = \sum_{k=0}^\infty \hat u_k \, L_k$ and $u' = \sum_{k=0}^\infty \hat u\sp{(1)}_k \, L_k$ are the orthogonal expansions of $u \in \HH^1(-1,1)$ and its weak derivative with respect to the basis $(L_k)_{k=0}^\infty$ of Legendre polynomials.
See \cite[Eq.~(2.4.22)]{CHQZ-I}---the first plus sign there is a typo---, \cite[Eq.~(2.13)]{Guo:2000a} and \cite[Lem.~3.4]{Figueroa:arXiv2015} for spectral differentiation formulas for Chebyshev, Gegenbauer and Zernike orthogonal polynomial expansions.
Whereas in one and two dimensions these particular bases of orthogonal polynomials are known to satisfy a wealth of simple identities so as to make spectral differentiation formulas simple to derive, that might not be the case for known explicit orthogonal polynomial bases $\LL^2_\alpha$ with $d \geq 3$ (cf.\ the example bases in \cite[\S~5.2]{DunklXu:2014}).

In this work we introduce a streamlined technique to prove \autoref{thm:lossy} which circumvents the need for spectral differentiation formulas and actually dispenses with the usage of bases of orthogonal polynomials altogether, focusing instead on orthogonal polynomial spaces; that is, spaces of polynomials of a certain degree orthogonal to all polynomials of lower degree (cf.\ \eqref{OPS} and the opening remarks of \cite[Ch.~3]{DunklXu:2014}).
In this way we can settle our main result seamlessly for any dimension.

The outline of this article is as follows.
In \autoref{sec:OP-WSS} we introduce some necessary additional notation, orthogonal polynomial spaces and some known properties of their members and their associated projectors.
The core of this work lies in \autoref{sec:id-diff}, in which we prove preliminary results concerning orthogonal polynomial spaces and their projectors.
Finally, in \autoref{sec:main} we bound a differentiation-projection commutator, prove our main result \autoref{thm:lossy} and an interpolation corollary and wrap up with some general remarks and a brief conclusion.

We finish this introductory section noting that have we omitted the dimension $d$ from the notation of $W_\alpha$, $\LL^2_\alpha$, etc.\ and will mostly continue to do so in order to avoid cluttering and because all of our arguments will be dimension-independent.

\section{Orthogonal polynomials and weighted Sobolev spaces}\label{sec:OP-WSS}

We denote by $\natural$ the set of strictly positive integers and $\natural_0 := \{0\} \cup \natural$.
Members of $[\natural_0]^d$ will be called multi-indices and for every multi-index $\gamma \in [\natural_0]^d$, point $x \in \Real^d$ and (strongly or weakly) differentiable enough complex-valued function $f$ defined on some open set of $\Real^d$ we shall write $\abs{\gamma} = \sum_{i=1}^d \gamma_i$, $x^\gamma = \prod_{i=1}^d x_i^{\gamma_i}$ and $\partial_\gamma f = \partial^{\abs{\gamma}} f/(\partial x_1^{\gamma_1} \dotsm \partial x_d^{\gamma_d})$.
We will denote by $\abs{\cdot}_{\HH^k_\alpha}$ the seminorm defined as the square root of $u \mapsto \norm{\nabla_k u}_\alpha^2 = \sum_{\abs{\gamma} = k} \binom{k}{\gamma} \norm{\partial_\gamma u}_\alpha^2$, where $\binom{k}{\gamma} = k!/(\gamma_1! \dotsm \gamma_d!)$ is the number of times the multi-index $\gamma$ of order $k$ appears in the $k$-dimensional array-valued $\nabla_k u$.
This seminorm is of course equivalent to the common choice in which the $\binom{k}{\gamma}$ are all replaced by $1$ yet better suits some induction arguments on the order of differentiation we make below.

Let $\mathcal{V}\sp{\alpha}_k$ be the space of orthogonal polynomials of degree $k$ with respect to the weight $W_\alpha$ (cf.\ \cite[Def.~3.1.1]{DunklXu:2014}); i.e.,
\begin{equation}\label{OPS}
\mathcal{V}\sp{\alpha}_k := \left\{ p \in \poly^d_k \mid (\forall\,q \in \poly^d_{k-1})\ \langle p, q\rangle_\alpha = 0 \right\}.
\end{equation}
If $k < 0$ we adopt the convention $\poly^d_k = \{0\}$ and so $\mathcal{V}\sp{\alpha}_k = \{0\}$.
As $W_\alpha$ is centrally symmetric, it transpires from \cite[Th.~3.3.11]{DunklXu:2014} that for all $k \in \natural_0 = \{0, 1, 2, \dotsc\}$ there holds the following parity relation:
\begin{equation}\label{parity}
(\forall \, p_k \in \mathcal{V}\sp{\alpha}_k)\ (\forall \, x \in B^d) \quad p_k(-x) = (-1)^k p_k(x).
\end{equation}
Let $\proj\sp{\alpha}_k$ denote the orthogonal projection from $\LL^2_\alpha$ onto $\mathcal{V}\sp{\alpha}_k$.
From \cite[Th.~3.2.18]{DunklXu:2014}, $\poly^d_n = \bigoplus_{k=0}^n \mathcal{V}\sp{\alpha}_k$ and $\LL^2_\alpha = \bigoplus_{k=0}^\infty \mathcal{V}\sp{\alpha}_k$, whence
\begin{equation}\label{spanning-consequence}
(\forall\,n\in\natural_0) \quad S\sp{\alpha}_n = \sum_{k=0}^n \proj\sp{\alpha}_k
\qquad\text{and}\qquad
(\forall\,u\in \LL^2_\alpha) \quad u = \sum_{k=0}^\infty \proj\sp{\alpha}_k(u).
\end{equation}
We will denote the entrywise application of $S\sp{\alpha}_n$ to $\LL^2_\alpha$-valued vectors and higher-order tensors by $S\sp{\alpha}_n$ as well (cf.\ \autoref{cor:diffProjComm-r} below).

From \cite[Eq.~(5.2.3) and Th.~8.1.3]{DunklXu:2014} and straightforward algebraic manipulation it is readily computed that the members of $\mathcal{V}\sp{\alpha}_k$ are eigenfunctions of the second order differential operator $p \mapsto -W_\alpha^{-1} \div\left( W_{\alpha+1} \nabla p \right) - \sum_{1 \leq i<j \leq d} D_{i,j}^2 p$, where $D_{i,j}$ denotes the first order angular differential operator $x_i \partial_j - x_j \partial_i$ \cite[\S~1.8]{DaiXu:2013}, with associated eigenvalue $k(k+d+2\alpha)$.
By integration by parts the following integral form follows:
\begin{multline}\label{weak-EV}
(\forall\,p_k \in \mathcal{V}\sp{\alpha}_k)\ \left(\forall\,q\in\CC^1(\overline{B^d})\right)\\
\langle \nabla p_k, \nabla q \rangle_{\alpha+1} + \sum_{1 \leq i < j \leq d} \left\langle D_{i,j} p_k, D_{i,j} q \right\rangle_\alpha
= k(k+d+2\alpha) \langle p_k, q \rangle_\alpha.
\end{multline}

\begin{remark}\label{rem:easy}
Together with appropriate density results, \eqref{weak-EV} implies that a member of $\mathcal{V}\sp{\alpha}_k$ is automatically also a member of an orthogonal polynomial subspace with respect to a Sobolev-type inner product involving the weaker weight $W_{\alpha+1}$ to control the gradient and, if $d \geq 2$, additional control for the angular derivatives.
In the $d = 1$ case, measuring the projection error in this induced non-uniformly weighted Sobolev space and its generalizations to higher degree of weak differentiation turns out to follow the trigonometric case much more closely (cf.\ \cite[Th.~2.1]{GW:2004} in the one-dimensional case with not necessarily symmetric Jacobi weights).
\end{remark}

\begin{lemma}\label{lem:density}
Let $d \in \natural$, $\alpha > -1$ and $m \in \natural_0$.
Then, $\CC^\infty(\overline{B^d})$ is dense in $\HH^m_\alpha$.
\begin{proof}
This follows from \cite[Rem.~11.12.(iii)]{Kufner:1985} upon the realization that $W_\alpha$ is bounded from above and below by positive multiples of $\operatorname{dist}(\cdot,\partial B^d)$.
\end{proof}
\end{lemma}

We cite from \cite[Cor.~2.7 and Lem.~2.11]{Figueroa:arXiv2015} the following $\LL^2_\alpha$ bound on the $S\sp{\alpha}_n$ projection error and an inverse or Markov-type inequality:

\begin{lemma}\label{lem:L2-projection-error}
For all $\alpha > -1$, $d \in \natural$ and $l \in \natural_0$ there exists a positive constant $C = C(\alpha,d,l)$ such that
\begin{equation*}
(\forall\,n\in\natural_0) \ (\forall\,u\in\HH^l_\alpha) \quad
\norm{u - S\sp{\alpha}_n(u)}_\alpha \leq C (n+1)^{-l} \norm{u}_{\HH^l_\alpha}.
\end{equation*}
\end{lemma}

\begin{lemma}\label{lem:Markov}
For $\alpha > -1$ and $d \in \natural$ there exists a positive constant $C = C(\alpha,d) > 0$ such that
\begin{equation*}
(\forall\,n\in\natural_0) \ (\forall\,p_n\in\poly^d_n) \quad
\norm{\nabla p_n}_\alpha \leq C n^2 \norm{p_n}_\alpha.
\end{equation*}
\end{lemma}

\section{Connections between orthogonal polynomials spaces and their projectors}\label{sec:id-diff}

The following proposition collects results concerning relations between spaces of orthogonal polynomials and their associated projectors not involving differentiation.

\begin{proposition}\label{pro:id-shift}
Let $\alpha > -1$ and $d \in \natural$.
\begin{enumerate}
\item\label{it:weighted-id-downshift} Let $p_k \in \mathcal{V}\sp{\alpha+1}_k$.
Then, $(1-\norm{\cdot}^2) p_k \in \mathcal{V}\sp{\alpha}_k \oplus \mathcal{V}\sp{\alpha}_{k+2}$.
\item\label{it:unnamed} Let $q_k \in \mathcal{V}\sp{\alpha}_k$.
Then, $q_k = \proj\sp{\alpha+1}_{k-2}(q_k) + \proj\sp{\alpha+1}_k(q_k)$.
\item\label{it:proto-id-shift} Let $u \in \LL^2_\alpha$.
Then, $\proj\sp{\alpha+1}_k(u) = \proj\sp{\alpha+1}_k\left( \proj\sp{\alpha}_k(u) + \proj\sp{\alpha}_{k+2}(u) \right)$.
\item\label{it:id-shift} Let $u \in \LL^2_\alpha$.
Then,
\begin{equation*}
\proj\sp{\alpha+1}_k(u)
= \proj\sp{\alpha}_k(u) + \proj\sp{\alpha+1}_k \circ \proj\sp{\alpha}_{k+2}(u) - \proj\sp{\alpha+1}_{k-2} \circ \proj\sp{\alpha}_k(u).
\end{equation*}
\end{enumerate}
\begin{proof}
Given $q \in \poly^d_{k-1}$, $\langle (1-\norm{\cdot}^2) p_k, q \rangle_\alpha = \langle p_k, q \rangle_{\alpha+1} = 0$ by definition \eqref{OPS}.
Also, by the parity relation \eqref{parity}, $(1-\norm{\cdot}^2) p_k \perp_\alpha \mathcal{V}\sp{\alpha}_{k+1}$.
Therefore part \ref{it:weighted-id-downshift} stems from \eqref{spanning-consequence}.
An analogous argument accounts for part \ref{it:unnamed}.
Part \ref{it:proto-id-shift} comes from the fact that given $p_k \in \mathcal{V}\sp{\alpha+1}_k$,
\begin{multline*}
\langle \proj\sp{\alpha+1}_k(u), p_k \rangle_{\alpha+1}
= \langle u, p_k \rangle_{\alpha+1}
= \langle u, (1-\norm{\cdot}^2) p_k \rangle_\alpha\\
\stackrel{\text{\ref{it:weighted-id-downshift}}}{=} \langle \proj\sp{\alpha}_k(u) + \proj\sp{\alpha}_{k+2}(u), (1-\norm{\cdot}^2) p_k \rangle_\alpha
= \langle \proj\sp{\alpha}_k(u) + \proj\sp{\alpha}_{k+2}(u), p_k \rangle_{\alpha+1}.
\end{multline*}
Part \ref{it:id-shift} is obtained from adding and subtracting $\proj\sp{\alpha+1}_{k-2}(\proj\sp{\alpha}_k(u))$ to the right hand side of part \ref{it:proto-id-shift} and using part \ref{it:unnamed}.
\end{proof}
\end{proposition}

We will now present another collection of results, this time involving differentiation.
To this end we introduce the first order differentiation operator $d\sp{\alpha}_j$, $\alpha > -1$ and $j \in \{1, \dotsc, d\}$, by
\begin{equation*}
d\sp{\alpha}_j q(x)
:= -W_\alpha(x)^{-1} \frac{\partial}{\partial x_j} \left( W_{\alpha+1}(x) \, q(x) \right)
= -(1-\norm{x}^2) \, \partial_j q(x) + 2 (\alpha+1) \, x_j \, q(x).
\end{equation*}

\begin{proposition}\label{pro:diff-shift}
Let $\alpha > -1$, $d \in \natural$ and $j \in \{1, \dotsc, d\}$.
\begin{enumerate}
\item\label{it:DARD} $d\sp{\alpha}_j$ maps $\poly^d_k$ into $\poly^d_{k+1}$.
\item\label{it:differentiationAdjoint} Given $p, q \in \CC^1(\overline{B^d})$, $\langle \partial_j p, q \rangle_{\alpha+1} = \langle p, d\sp{\alpha}_j q \rangle_\alpha$.
\item\label{it:DALP} Let $r_k \in \mathcal{V}\sp{\alpha+1}_k$.
Then, $d\sp{\alpha}_j(r_k) \in \mathcal{V}\sp{\alpha}_{k+1}$.
\item\label{it:proto-diff-shift} Let $p_k \in \mathcal{V}\sp{\alpha}_k$.
Then, $\partial_j p_k \in \mathcal{V}\sp{\alpha+1}_{k-1}$.
\item\label{it:diff-shift} Let $u \in \CC^1(\overline{B^d})$.
Then, $\partial_j\proj\sp{\alpha}_k(u) = \proj\sp{\alpha+1}_{k-1}(\partial_j u)$.
\end{enumerate}
\begin{proof}
Part \ref{it:DARD} is straightforward.
Part \ref{it:differentiationAdjoint} is obtained by integration by parts and noticing that no boundary term appears on account of $(1-\norm{\cdot}^2)^{\alpha+1}$ vanishing on the boundary of $B^d$.

Given $r_k \in \mathcal{V}\sp{\alpha+1}_k$, by part \ref{it:DARD}, $d\sp{\alpha}_j(r_k) \in \poly^d_{k+1}$, and, on account of part \ref{it:differentiationAdjoint}, it is $\LL^2_\alpha$-orthogonal to $\poly^d_k$, whence part \ref{it:DALP}.
An analogous argument accounts for part \ref{it:proto-diff-shift}.

Given $u \in \CC^1(\overline{B^d})$, by part \ref{it:proto-diff-shift}, $\partial_j \proj\sp{\alpha}_k(u) \in \mathcal{V}\sp{\alpha+1}_{k-1}$.
Part \ref{it:diff-shift} then comes about from the fact that for all $r \in \mathcal{V}\sp{\alpha+1}_{k-1}$,
\begin{equation*}
\langle \partial_j\proj\sp{\alpha}_k(u), r \rangle_{\alpha+1}
\stackrel{\text{\ref{it:differentiationAdjoint}}}{=} \langle \proj\sp{\alpha}_k(u), d\sp{\alpha}_j r \rangle_\alpha
\stackrel{\text{\ref{it:DALP}}}{=} \langle u, d\sp{\alpha}_j r \rangle_\alpha
\stackrel{\text{\ref{it:differentiationAdjoint}}}{=} \langle \partial_j u, r \rangle_{\alpha+1}.
\end{equation*}
\end{proof}
\end{proposition}

\begin{remark}[Shift operators]\label{rem:shift}
Part \ref{it:DALP} of \autoref{pro:diff-shift} means that $d\sp{\alpha}_j$ is a \emph{backward shift/degree raising} operator in the sense of \cite{KLS:2010}.
Similarly, by part \ref{it:proto-diff-shift}, $\partial_j$ is a \emph{forward shift/degree lowering} operator (see also \eqref{Jacobi-diff-shift} below).
\end{remark}

Inasmuch as it allows for quantifying a ``wrong'' ($\LL^2_\alpha$) norm of a member of a space of orthogonal polynomials ($\mathcal{V}\sp{\alpha+1}_k$), the following result is distantly related to \cite[Eq.~(4.43)]{Figueroa} and \cite[Prop.~3.12]{Figueroa:arXiv2015} in the $d =1$ and $d = 2$ cases, respectively.

\begin{proposition}\label{pro:extendedOrthogonality}
Let $\alpha > -1$, $d \in \natural$ and $k \in \natural_0$.
Then, for all $p, q \in \mathcal{V}\sp{\alpha+1}_k$,
\begin{equation*}
\langle p, q \rangle_\alpha = \left(\frac{k+d/2}{\alpha+1} + 1\right)\langle p, q \rangle_{\alpha+1}.
\end{equation*}
\begin{proof}
We start with the observation that if $s$ is a homogeneous polynomial of degree $k$---that is, of the form $s(x) = \sum_{\abs{\gamma}=n} c_\gamma x^\gamma$---, it satisfies $x \cdot \nabla s(x) = k\, s(x)$, which also goes on to show that the $x \cdot \nabla$ operator exactly preserves the degree of any $d$-variate polynomial.

Let $p, q \in \mathcal{V}\sp{\alpha+1}_k$.
As every member of $\mathcal{V}\sp{\alpha+1}_k$ is a linear combination of homogeneous polynomials of degree ranging from $0$ to $k$, there exists a homogeneous polynomial $s_p$ of degree $k$ such that $p - s_p \in \poly^d_{k-1}$ and hence $x \cdot \nabla p - x \cdot \nabla s_p \in \poly^d_{k-1}$.
Thus,
\begin{equation}\label{origin-of-k}
\langle x \cdot \nabla p, q\rangle_{\alpha+1}
= \langle x \cdot \nabla s_p, q \rangle_{\alpha+1}
= k \langle s_p, q \rangle_{\alpha+1}
= k \langle p, q \rangle_{\alpha+1}.
\end{equation}
Using the facts that $\nabla (1-\norm{x}^2)^{\alpha+1} = -2(\alpha+1) (1-\norm{x}^2)^\alpha x$, $\operatorname{div}(x) = d$, integration by parts and \eqref{origin-of-k}, which of course is still valid if the roles of $p$ and $q$ are interchanged,
\begin{multline*}
2(\alpha+1)\int_{B^d} p(x) \overline{q(x)} \norm{x}^2 (1-\norm{x}^2)^\alpha \dd x
= \int_{B^d} \mathrm{div}\left(p(x) \overline{q(x)} x\right) (1-\norm{x}^2)^{\alpha+1} \dd x\\
= \left( \langle x \cdot \nabla p, q\rangle_{\alpha+1} + \langle p, x \cdot \nabla q \rangle_{\alpha+1} + d \langle p, q \rangle_{\alpha+1} \right)
= (2k+d) \langle p, q \rangle_{\alpha+1}.
\end{multline*}
The desired result then follows from the fact that $(1-\norm{x}^2)^\alpha = \norm{x}^2 (1-\norm{x}^2)^\alpha + (1-\norm{x}^2)^{\alpha+1}$.
\end{proof}
\end{proposition}

\begin{remark}[Relations with identities satisfied by bases]\label{rem:variants-with-bases}
In the one-dimensional case ($d = 1$), $\mathcal{V}\sp{\alpha}_k = \operatorname{span}(\{P\sp{(\alpha,\alpha)}_k\})$, where the $P\sp{(\alpha,\alpha)}_k$ are Jacobi polynomials \cite[Ch.~4]{Szego:1975}.
Then, from the ``id-shift'' identity (a combination of (6.4.21) and (6.4.23) of \cite{AAR:1999}; it must be slightly modified if $\alpha = -1/2$ and $k = 0$)
\begin{equation}\label{Jacobi-id-shift}
P\sp{(\alpha,\alpha)}_k = \frac{(k+2\alpha+1)(k+2\alpha+2)}{(2k+2\alpha+1)(2k+2\alpha+2)} P\sp{(\alpha+1,\alpha+1)}_k - \frac{k+\alpha}{2(2k+2\alpha+1)} P\sp{(\alpha+1,\alpha+1)}_{k-2},
\end{equation}
it is possible to furnish alternative proofs of parts \ref{it:unnamed} and \ref{it:proto-id-shift} of \autoref{pro:id-shift} and hence of its part \ref{it:id-shift}.
In that rough sense \autoref{pro:id-shift} corresponds to \eqref{Jacobi-id-shift}.
Similarly \cite[Eq.~(4.21.7)]{Szego:1975},
\begin{equation}\label{Jacobi-diff-shift}
{P\sp{(\alpha,\alpha)}_k}' = \frac{k+2\alpha+1}{2} P\sp{(\alpha+1,\alpha+1)}_{k-1},
\end{equation}
allows for proving part \ref{it:diff-shift} of \autoref{pro:diff-shift} and so, again in a rough sense, \autoref{pro:diff-shift} corresponds to \eqref{Jacobi-diff-shift}.
Using \eqref{Jacobi-id-shift} and explicit formulas for the norms of Jacobi polynomials (cf.\ \cite[Eq.~(4.3.3)]{Szego:1975}) it is possible to reconstruct \autoref{pro:extendedOrthogonality}, although the necessary computations are not short.

In the two-dimensional case, $\mathcal{V}\sp{\alpha}_k = \operatorname{span}(\{ P\sp{(\alpha)}_{m,n} \mid m+n = k \})$, where each $P\sp{(\alpha)}_{m,n}$ is a Zernike polynomial \cite{Wunsche:2005}.
Then, the identities \eqref{Jacobi-id-shift} and \eqref{Jacobi-diff-shift} find appropriate analogues in \cite[Eq.~(3.12)]{Figueroa:arXiv2015} and \cite[Eq.~(5.3)]{Wunsche:2005}, respectively.
\end{remark}

\section{Proof of the main result and an interpolation corollary}\label{sec:main}

We can now bound a differentiation-projection commutator.

\begin{lemma}\label{lem:diffProjComm-1}
Let $\alpha > -1$, $d \in \natural$ and $l \in \natural$.
Then, there exists $C = C(\alpha,d,l) > 0$ such that for all $u \in \HH^l_\alpha$, $n \in \natural_0$ and $j \in \{1, \dotsc, d\}$,
\begin{equation*}
\norm{\partial_j S\sp{\alpha}_n(u) - S\sp{\alpha}_n(\partial_j u)}_\alpha
\leq C (n+1)^{3/2-l} \norm{\partial_j u}_{\HH^{l-1}_\alpha}.
\end{equation*}
\begin{proof}
Let us first assume that $u \in \CC^\infty(\overline{B^d})$.
Combining part \ref{it:id-shift} of \autoref{pro:id-shift} and part \ref{it:diff-shift} of \autoref{pro:diff-shift}, we obtain
\begin{equation}\label{id-diff}
\partial_j\proj\sp{\alpha}_{k+1}(u) - \proj\sp{\alpha}_k(\partial_j u)
= \proj\sp{\alpha+1}_k \circ \proj\sp{\alpha}_{k+2}(\partial_j u) - \proj\sp{\alpha+1}_{k-2} \circ \proj\sp{\alpha}_k(\partial_j u).
\end{equation}
Using \eqref{spanning-consequence} to express $S\sp{\alpha}_n$ in terms of the $\proj\sp{\alpha}_k$, using \eqref{id-diff}, noticing that a telescoping sum results and using part \ref{it:unnamed} of \autoref{pro:id-shift} to expand an appearance of $\proj\sp{\alpha}_n(\partial_j u) \in \mathcal{V}\sp{\alpha}_n$,
\begin{multline}\label{commutator}
\partial_j S\sp{\alpha}_n(u) - S\sp{\alpha}_n(\partial_j u)
= \sum_{k=0}^n \partial_j \proj\sp{\alpha}_k(u) - \sum_{k=0}^n \proj\sp{\alpha}_k(\partial_j u)\\
= \sum_{k=0}^{n-1} \left( \partial_j \proj\sp{\alpha}_{k+1}(u) - \proj\sp{\alpha}_k(\partial_j u) \right) - \proj\sp{\alpha}_n(\partial_j u)\\
= \proj\sp{\alpha+1}_{n-2} \circ \proj\sp{\alpha}_n(\partial_j u) + \proj\sp{\alpha+1}_{n-1} \circ \proj\sp{\alpha}_{n+1}(\partial_j u) - \proj\sp{\alpha}_n(\partial_j u)\\
= \proj\sp{\alpha+1}_{n-1} \circ \proj\sp{\alpha}_{n+1}(\partial_j u) - \proj\sp{\alpha+1}_n \circ \proj\sp{\alpha}_n(\partial_j u).
\end{multline}
Now, by \autoref{pro:extendedOrthogonality}, the fact that $\norm[n]{\proj\sp{\alpha+1}_{n-1}}_{\mathcal{L}(\LL^2_{\alpha+1})} \leq 1$ and the fact that $\norm{\cdot}_{\alpha+1} \leq \norm{\cdot}_\alpha$ in $\LL^2_\alpha$ (because $W_{\alpha+1} \leq W_\alpha$) we have that for all $n \geq 1$,
\begin{equation}\label{good-enough-1}
\norm{\proj\sp{\alpha+1}_{n-1} \circ \proj\sp{\alpha}_{n+1}(\partial_j u)}_\alpha^2
\leq \frac{n+d/2+\alpha}{\alpha+1} \norm{\proj\sp{\alpha}_{n+1}(\partial_j u)}_\alpha^2.
\end{equation}
Of course, if $n = 0$, our conventions imply that $\norm{\proj\sp{\alpha+1}_{n-1} \circ \proj\sp{\alpha}_{n+1}(\partial_j u)}_\alpha^2 = 0$.
Analogous arguments show that for all $n \in \natural_0$,
\begin{equation}\label{good-enough-2}
\norm{\proj\sp{\alpha+1}_n \circ \proj\sp{\alpha}_n(\partial_j u)}_\alpha^2
\leq \frac{n+1+d/2+\alpha}{\alpha+1} \norm{\proj\sp{\alpha}_n(\partial_j u)}_\alpha^2.
\end{equation}
Taking the squared $\LL^2_\alpha$ norm of both ends of \eqref{commutator}, exploiting the $\LL^2_\alpha$ orthogonality of $\mathcal{V}\sp{\alpha+1}_{n-1}$ and $\mathcal{V}\sp{\alpha+1}_n$ (a consequence of the parity relation \eqref{parity}) and the bounds \eqref{good-enough-1} and \eqref{good-enough-2} we observe that
\begin{equation*}
\norm{\partial_j S\sp{\alpha}_n(u) - S\sp{\alpha}_n(\partial_j u)}_\alpha^2
\leq \frac{n+1+d/2+\alpha}{\alpha+1} \norm{\partial_j u - S\sp{\alpha}_{n+2}(\partial_j u)}_\alpha^2.
\end{equation*}
As $\partial_j u \in \HH^{l-1}_\alpha$, we can appeal to \autoref{lem:L2-projection-error} to obtain the desired result for $u \in \CC^\infty(\overline{B^d})$ after realizing that there exists a constant $\tilde C$ depending only on $\alpha$, $d$ and $l$ such that $\frac{n+1+d/2+\alpha}{\alpha+1} ((n+3)^{-(l-1)})^2 \leq \tilde C (n+1)^{3-2l}$ for all $n \in \natural_0$.
The general result then follows via the density result in \autoref{lem:density}.
\end{proof}
\end{lemma}

\begin{corollary}\label{cor:diffProjComm-r}
Let $\alpha > -1$, $d \in \natural$ and $r, l \in \natural$ with $r \leq l$.
Then, there exists $C = C(\alpha,d,l,r) > 0$ such that for all $u \in \HH^l_\alpha$ and $n \in \natural_0$,
\begin{equation*}
\norm{\nabla_r S\sp{\alpha}_n(u) - S\sp{\alpha}_n(\nabla_r u)}_\alpha \leq C (n+1)^{2r-1/2-l} \norm{u}_{\HH^l_\alpha}.
\end{equation*}
\begin{proof}
Let us first note that iterating \autoref{lem:Markov} we find that for all $r \in \natural$ there exists $C > 0$ depending on $\alpha$, $d$ and $r$ such that
\begin{equation}\label{iteratedMarkov}
(\forall\,n\in\natural_0)\ (\forall\,p\in\poly^d_n) \quad \abs{p}_{\HH^r_\alpha} \leq C n^{2 r} \norm{p}_\alpha.
\end{equation}

We will now operate by induction on $r$.
Taking the square root of the sum with respect to $j$ of the square of both sides of the inequality in \autoref{lem:diffProjComm-1} the case $r = 1$ follows almost immediately.
Let us suppose now that our desired result holds for some $r \in \{1, \dotsc, l\}$ and that $r+1 \leq l$.
Then, for all $j \in \{1, \dotsc, d\}$, by the triangle inequality,
\begin{equation*}
\norm{\nabla_r \partial_j S\sp{\alpha}_n(u) - S\sp{\alpha}_n(\nabla_r \partial_j u)}_\alpha
\leq \abs{\partial_j S\sp{\alpha}_n(u) - S\sp{\alpha}_n(\partial_j u)}_{\HH^r_\alpha} + \norm{\nabla_r S\sp{\alpha}_n(\partial_j u) - S\sp{\alpha}_n(\nabla_r \partial_j u)}_\alpha.
\end{equation*}
By \eqref{iteratedMarkov} and \autoref{lem:diffProjComm-1} the first term is bounded by an appropriate constant times $n^{2r} (n+1)^{3/2-l} \norm{\partial_j u}_{\HH^{l-1}_\alpha}$.
By the induction hypothesis and the fact that $\partial_j u \in \HH^{l-1}_\alpha$ the second term is bounded by an appropriate constant times $(n+1)^{2r-1/2-(l-1)} \norm{\partial_j u}_{\HH^{l-1}_\alpha}$.
Then the desired result in the $r+1$ case follows from summing up with respect to $j$ and standard inequalities connecting vector $1$- and $2$-norms.
\end{proof}
\end{corollary}

We are now in a position prove our main result, \autoref{thm:lossy}, and the interpolation \autoref{cor:interpolatedLossy}.
As those proofs are almost completely analogous to those of Theorem~3.9 and Corollary~3.10 of \cite{Figueroa:arXiv2015} we only sketch them here.

\begin{proof}[Proof of \autoref{thm:lossy}] For every $k \in \{1, \dotsc, r\}$,
\begin{equation*}
\abs{u - S\sp{\alpha}_N(u)}_{\HH^k_\alpha}^2
\leq 2 \norm{\nabla_k u - S\sp{\alpha}_N(\nabla_k u)}_\alpha^2 + 2 \norm{S\sp{\alpha}_N(\nabla_k u) - \nabla_k S\sp{\alpha}_N(u)}_\alpha^2.
\end{equation*}
We bound the first term using \autoref{lem:L2-projection-error} and the second using \autoref{cor:diffProjComm-r} and the desired result follows upon summing up with respect to $k$ and taking the square root.
\end{proof}

Given $m \in \natural_0$ and $\theta \in (0,1)$ we define $\HH^{m+\theta}_\alpha$ by complex interpolation \cite[\P7.51--52]{AF:2003}:
\begin{equation}\label{InterpolatedSobolev}
\HH^{m+\theta}_\alpha := \left[\HH^m_\alpha, \HH^{m+1}_\alpha\right]_\theta.
\end{equation}

\begin{corollary}\label{cor:interpolatedLossy}
Let $\alpha > -1$, $d \in \natural$ and $r, l \geq 0$ with $r \leq l$.
Then, there exists $C = C(\alpha,d,l,r) > 0$ such that for all $u \in \HH^l_\alpha$ and $n \in \natural_0$,
\begin{equation*}
\norm{u - S\sp{\alpha}_n(u)}_{\HH^r_\alpha} \leq C n^{e(l,r)} \norm{u}_{\HH^l_\alpha}
\quad\text{where}\quad
e(l,r) = \begin{cases}
3/2 \, r - l & \text{if } 0 \leq r \leq 1,\\
2\, r - 1/2 - l & \text{if } r \geq 1.
\end{cases}
\end{equation*}
\begin{proof}
The desired bound on the operator norm of $T\sp{\alpha}_{n,l,r} \colon \HH^l_\alpha \to \HH^r_\alpha$ defined by $T\sp{\alpha}_{n,l,r} := I - S\sp{\alpha}_n$ (with $I$ being the identity operator) holds when $r$ and $l$ are integers from \autoref{lem:L2-projection-error} in the $r = 0$ case and \autoref{thm:lossy} in the $r \in \natural$ case.
The non-integer cases then follow by using the exact interpolation and reiteration theorems.
\end{proof}
\end{corollary}

\begin{remark}[Real interpolation]
Just as it was remarked upon in the $d = 2$ case in \cite{Figueroa:arXiv2015}, essentially the same argument used in \autoref{cor:interpolatedLossy} would work if we used real instead of complex interpolation to define the weighted Sobolev spaces with non-integer differentiation parameter in \eqref{InterpolatedSobolev}.
\end{remark}

\begin{remark}[On the optimality of the main result]
There are four parameters in our main result, \autoref{thm:lossy}: The dimension $d \in \natural$, the weight parameter $\alpha \in (-1, \infty)$, the regularity parameter of the function being approximated $l \in \natural$ and the regularity parameter of the norm measuring the residual $r \in \{1, \dotsc, l\}$.
We will say that \autoref{thm:lossy} is optimal if the power on $N$ in \eqref{lossyPreview} cannot be lowered.
We are aware of optimality proofs in the cases $(d,\alpha,l,r) = (1,-1/2,1,1)$ \cite[pp.~76, 78]{CQ:1982}, $(d,\alpha,l,r) = (1,0,1,1)$ \cite[p.~285]{CHQZ-I}, $(d,\alpha,l,r) \in \{2\} \times (-1,\infty) \times \natural \times \{1\}$ \cite[Th.~3.13]{Figueroa:arXiv2015} (the latter can be adapted to $(d,\alpha,l,r) \in \{1\} \times (-1,\infty) \times \natural \times \{1\}$).
All those proofs exploit a number of simple identities satisfied by particular bases of orthogonal polynomials.
Notice also that all those parameter regimes have $r = 1$, arguably the most important $r$ in \autoref{thm:lossy} because of its connection with the analysis of weak forms of second order PDE.
In \cite{Figueroa:arXiv2015} numerical experiments were used to support the conjecture that \autoref{thm:lossy} is also optimal for $(d,\alpha,l,r) \in \{2\} \times (-1,\infty) \times \{(l,r) \in \natural \times \natural \mid r \leq l\}$.
For general $d$ we do not know of bases of $\mathcal{V}\sp{\alpha}_k$ satisfying identities (particularly regarding differentiation) simple enough so as to enable us to completely extend the optimality proofs mentioned above.
Nevertheless, always in the $r = 1$ case, we managed to generalize the techniques used in \cite{Figueroa:arXiv2015} for $(\alpha,l)$ in a certain proper subset of its natural range $(-1,\infty) \times \natural$.
The arguments behind this partial result being rather involved, depending on explicit identities satisfied by Jacobi polynomials and thus out of character with the rest of this work, we decided against including them here.
\end{remark}

\subsection*{Conclusion}

We have proved our desired ``lossy'' (as compared to the unweighted trigonometric case) bound \autoref{thm:lossy} and did so without recourse to special identities satisfied by particular bases of orthogonal polynomials, arguing instead in terms of orthogonal polynomial spaces.

We certainly expect the main sequence of results in \autoref{sec:id-diff} and \autoref{sec:main} to extend to a wider class of of reflection invariant weights.
If we focused on Gegenbauer-type weights it was mostly on account of their importance in applications and the ready availability of \autoref{lem:density}, \autoref{lem:L2-projection-error} and \autoref{lem:Markov}.

\bibliographystyle{abbrvnat}
\small
\bibliography{bpa-refs}
\normalsize

\bigskip
\bigskip

\end{document}